\documentclass{amsart}

\usepackage[all]{xy}
\usepackage{xcolor}

\newtheorem{theorem}{Theorem}[section]

\theoremstyle{definition}
\newtheorem{definition}[theorem]{Definition}
\newtheorem{example}[theorem]{Example}
\newtheorem{proposition}[theorem]{Proposition}

\theoremstyle{remark}

\theoremstyle{notation}

\numberwithin{equation}{section}

\begin{document}
\author{Yanlong Hao}
\address{School of Mathematical Sciences,
         Nankai University,
         Tianjin 300071, P.R.China}
\email{haoyanlong13@mail.nankai.edu.cn}
\author{Xiugui Liu*}
\address{School of Mathematical Sciences and LPMC,
         Nankai University,
         Tianjin 300071, P.R.China}
\email{xgliu@nankai.edu.cn}
\author{Qianwen Sun}
\address{School of Mathematical Sciences,
         Nankai University,
         Tianjin 300071, P.R.China}
\email{qwsun13@mail.nankai.edu.cn}
\thanks{*The second author was supported in part by the National Natural Science Foundation of China (No. 11571186), and the Scientific Research Foundation for the Returned Overseas Chinese Scholars, State Education Ministry.}
\title{Two-stage spaces and the torus rank conjecture}
\maketitle
$\textbf{Abstract}$: In this note, we give some new families of two-stage spaces for which the torus rank conjecture is affirmed.

$\textbf{Key words}$: Rational homotopy theory, two-stage Sullivan algebra, torus rank conjecture.
\section{Introduction}
The torus rank conjecture is one of the long-standing problems in algebraic topology. The conjecture says that if a nilpotent finite dimensional CW-complex admits an almost free $T^n$ action, then 
$${\rm{dim}}~H(X;\mathbb{Q})\geq 2^n.$$

We shall henceforth assume that all spaces are 1-connected, finite cell complexes. Now, the conjecture is proved for $n\leq 3$ (\cite{C-V}). There are many families of spaces for which the conjecture holds. For instance it is proved that the conjecture holds for homogeneous spaces of a compact Lie group, and spaces satisfying the hard Lefschetz properties (\cite{A-C-P}).

There is an algebraic version of this conjecture \cite{H-S}, if $M$ is an $m$-dimensional compact manifold, then there is an almost free $T^n$ action on $M$ if there is a relative Sullivan algebra of the form
$$(\Lambda(x_1,x_2,\cdots,x_n),0)\rightarrow (\Lambda(x_1,x_2,\cdots,x_n)\otimes \Lambda V,D)\rightarrow (\Lambda V, d),$$
where $|x_i|=2$ for $i=1,2,\cdots,n$, $(\Lambda V,d)$ is a minimal model of $M$ and the cohomology groups $H^*(\Lambda(x_1,x_2,\cdots,x_n)\otimes \Lambda V,D)$ are finite dimensional. Moreover, if the above condition holds, then $T^n$ acts freely on a finite CW-complex $X$ which has the same homotopy type as $M$, and if $m-n\neq 0$ mod $4$, we can choose $X$ to be a compact manifold. In this case, the relative Sullivan algebra is a model of the associated Borel fibration of the torus action.

A Sullivan algebra is said to be a two-stage algebra if $(\Lambda V,d)\cong (\Lambda (U\oplus W), d)$ such that $dU=0$ and $dW\subset \Lambda U$. We call
the isomorphism $V\cong U\oplus W$ a two-stage decomposition of $(\Lambda V,d)$. If the minimal model of a space $X$ is a two-stage Sullivan algebra, we say that $X$ is a two-stage space.

A Sullivan algebra $(\Lambda W,d)$ with $W$
and $H(\Lambda W, d)$ both finite-dimensional are called elliptic, and a space $X$ with an
elliptic minimal model is called an elliptic space. Topologically, this means that
both $\pi(X)\otimes \mathbb{Q}$ and $H(X;\mathbb{Q})$ are finite dimensional.

In this paper, we will show the following theorem.
\begin{theorem}\label{1}
Let $(\Lambda V,d)$ be a two-stage algebra such that $V$ is a finite dimensional graded vector space. If $k_{V}\leq 3$, the torus rank conjecture holds for this algebra, where $k_{V}$ is defined in Section 2.
\end{theorem}
\section{Basics of rational homotopy theory}
At first, we recall some basics of rational homotopy theory, all of them can be found in \cite{Y-S}. We simply remark a few facts.

We recall that when $M$ is path connected, the Sullivan model of $M$ is a quasi-isomorphism:
$$m:(\Lambda V_{M},d)\rightarrow A_{PL}(M),$$
where $(\Lambda V_{M},d)$ is a Sullivan algebra.

A minimal model of a map $f:X\rightarrow Y$ is a quasi-isomorphism:
$$\varphi:(\Lambda V_{Y}\otimes \Lambda W,D)\rightarrow (\Lambda V_X,d),$$
where $(\Lambda V_Y,d)$ and $(\Lambda V_{X},d)$ are the models of $Y$ and $X$, respectively, and also $Dw\in (\Lambda^+V_Y\otimes \Lambda W)\oplus \Lambda^+W$ for $w\in W$.

Given a nilpotent fibration $F\hookrightarrow X\rightarrow Y$ where the involved spaces are all of finite type, then there is a  K-S extension model of the fibration given by
$$(\Lambda V_{Y},d)\hookrightarrow (\Lambda V_{Y}\otimes \Lambda V_{F},D)\rightarrow (\Lambda V_{F},d),$$
where $(\Lambda V_{Y}\otimes \Lambda V_{F},D)$ is a relative minimal model of the projection.

A pure Sullivan algebra is a special two-stage algebra $(\Lambda V,d)$, which has a two-stage decomposition $V=U\oplus W$ such that $U=V^{\rm{even}}$.

Let $(\Lambda V,d)$ be a two-stage algebra such that ${\rm{dim}}~V<\infty$ and ${\rm{dim}}~H(\Lambda V,d)<\infty.$
Let $V=U\oplus W$ be the two-stage decomposition of this algebra. We have that $W$ is concentrated in odd degrees. Let $d=d_{\tau}+\theta$, where $d_{\tau}v\in \Lambda V^{\rm{odd}}$ and $\theta v\in \Lambda^+ V^{{\rm{even}}}\cdot \Lambda V$ for $v\in V$. Then $(\Lambda V,d_{\tau})$ is also a two-stage algebra. We have a new two-stage decomposition $V=U'\oplus W'$  for the new algebra.

\begin{definition}\cite{B-G}
Suppose $(\Lambda V, d)$ is a two-stage minimal Sullivan algebra. We say that a two-stage decomposition $V=U\oplus W$
has maximal dimension, or that the two-stage decomposition displays $W$ with maximal
dimension, if, for any two-stage decomposition $V=\widetilde{U}\oplus \widetilde{W}$, we have ${\rm{dim}}~ W\geq {\rm{dim}}~\widetilde{W}.$
\end{definition}
\begin{definition}
$k_{V}={\rm{dim}}~U'$, where $V=U'\oplus W'$ is a two-stage decomposition of $(\Lambda V,d_{\tau})$, and has maximal dimension.
\end{definition}
\begin{example}
Let $(\Lambda(x_i,y_j,z,w),d)$ where $1\leq i,j\leq 6$, $|x_i|=3$, $|y_i|=2$, $dz=x_1x_2x_3x_4+x_5x_6y_1^3+y_2^6$ and $dw=y_3y_4y_5y_6$. We have $d_{\tau}z=x_1x_2x_3x_4$, hence $k_V=4$.
\end{example}
\begin{example}
It is also trivial that $k_V=0$ for a hyperelliptic Sullivan algebra. Here a hyperelliptic Sullivan algebra means a two-stage Sullivan algebra $(\Lambda V,d)$ such that $dV\in \Lambda^+V^{\rm{even}}\cdot \Lambda V$.
\end{example}
Then, we compute the homology dimension in a special case. The main result is then
established by reducing to this special case.
\begin{proposition}\label{p1}
Let $(\Lambda(U\oplus W ),d)$ be a two-stage, elliptic minimal Sullivan algebra with odd
degree generators only and suppose that $d: W\rightarrow \Lambda^2U$ is an isomorphism. Then
$${\rm{dim}}~H(\Lambda(U\oplus W),d)=\left\{
\begin{array}{ll}
2, & {\rm{dim}}~U=1,\\
6, & {\rm{dim}}~U=2,\\
36,& {\rm{dim}}~U=3.
\end{array}
\right.
$$
\end{proposition}
\begin{proof}
The desired result follows from a straightforward calculation.
\end{proof}
Let
$$\chi(k_{V})=2^{-k_{V}(k_{V}+1)/2}{\rm{dim}}~H(\Lambda (U\oplus W),D),$$
where $U=U^{\rm{odd}}$, ${\rm{dim}}~U=k_{V}$ and $D:W\rightarrow \Lambda^2 U$ is an isomorphism of vector spaces. Then $\chi(0)=1$, $\chi(1)=1$, $\chi(2)=3/4$ and $\chi(3)=9/16$.

As promised, the main result is now obtained by reducing the general two-stage cases to that of Proposition \ref{p1}.
\begin{theorem}\label{p3}
Suppose $(\Lambda V,d)$ is a two-stage, elliptic minimal Sullivan algebra. Then
$${\rm{dim}}~H(\Lambda(V),d)
\geq \chi(k_{V})2^{{\rm{dim}}~V^{\rm{odd}}-{\rm{dim}}~V^{\rm{even}}}.$$
\end{theorem}
\begin{proof}
At first, let $V=U\oplus W$ be the two-stage decomposition of $(\Lambda V,d)$ which has maximal dimension, then there is a K-S extension
$$(\Lambda V, d)\hookrightarrow (\Lambda (V\oplus \overline{U}), D)\rightarrow (\Lambda \overline{U}, 0),$$
where $D:\overline{U}\rightarrow U^{\rm{even}}$ is an vector space isomorphism.

This extension sequence has an associated Serre spectral sequence. This spectral sequence has $E_2$-term isomorphic to $H(\Lambda \overline{U}, 0)\otimes H(\Lambda V, d)$ and it converges to $H(\Lambda (V\oplus \overline{U}), D)$. Then we have
$${\rm{dim}}~ H(\Lambda \overline{U}, 0)\times {\rm{dim}}~ H(\Lambda V, d)\geq {\rm{dim}}~H(\Lambda (V\oplus \overline{U}), D).$$

We also have a K-S extension
$$(\Lambda (U\oplus\overline{U}),D)\hookrightarrow(\Lambda (V\oplus \overline{U}), D)\rightarrow (\Lambda (U^{\rm{odd}}\oplus W), d_{\tau}).$$ 
Since $H(\Lambda (U\oplus\overline{U}),D)=\mathbb{Q}$, we have $H(\Lambda (V\oplus \overline{U}), D)\simeq H(\Lambda (U^{\rm{odd}}\oplus W), d_{\tau})$.

Notice that
$${\rm{dim}}~H(\Lambda \overline{U}, 0)=2^{{\rm{dim}}~U^{\rm{even}}},$$
which implies
$${\rm{dim}}~H(\Lambda V, d)\geq 2^{-{\rm{dim}}~U^{\rm{even}}}{\rm{dim}}H(\Lambda (U^{\rm{odd}}\oplus W), d_{\tau}).$$

There is another K-S extension
$$(\Lambda (U^{\rm{odd}}\oplus W), d_{\tau})\hookrightarrow (\Lambda (U^{\rm{odd}}\oplus W\oplus \overline{W}), D)\rightarrow (\Lambda \overline{W}, 0),$$
where $D: \overline{W}\rightarrow \Lambda^2 U^{\rm{odd}}$ is an isomorphism of vector spaces. Notice that we have ${\rm{dim}}~H(\Lambda \overline{W}, 0)=2^{k_{V}(k_V-1)/2}$, and
$$(\Lambda (U^{\rm{odd}}\oplus W\oplus \overline{W}), D)\simeq (\Lambda (U^{\rm{odd}}\oplus \overline{W}),D)\otimes (\Lambda W,0).$$
Then we have
$${\rm{dim}}~H(\Lambda (U^{\rm{odd}}\oplus W), d_{\tau})\times {\rm{dim}}~H(\Lambda \overline{W}, 0)\geq {\rm{dim}}~H((\Lambda (U^{\rm{odd}}\oplus \overline{W}),D)\otimes (\Lambda W,0)),$$
which implies
$${\rm{dim}}~H(\Lambda (U^{\rm{odd}}\oplus W), d_{\tau})\geq 2^{{\rm{dim}} W-k_{V}(k_V-1)/2}{\rm{dim}}~H(\Lambda (U^{\rm{odd}}\oplus \overline{W}),D).$$
So we have
$${\rm{dim}}~H(\Lambda V, d)\geq 2^{{\rm{dim}}~W-k_{V}(k_V-1)/2-{\rm{dim}}~U^{\rm{even}}}{\rm{dim}}~H(\Lambda (U^{\rm{odd}}\oplus \overline{W}),D).$$

We also have ${\rm{dim}}~V^{\rm{odd}}-{\rm{dim}}~V^{\rm{even}}={\rm{dim}}~W+k_V-{\rm{dim}}~U^{\rm{even}}.$  Then
$${\rm{dim}}~H(\Lambda V, d)\geq 2^{{\rm{dim}}~V^{\rm{odd}}-{\rm{dim}}~V^{\rm{even}}}2^{-k_{V}(k_{V}+1)/2}{\rm{dim}}~H(\Lambda (U^{\rm{odd}}\oplus \overline{W}),D),$$
the desired result follows.
\end{proof}
Then we give an estimation of the torus rank of some two-stage minimal Sullivan algebras.
\begin{proposition}\label{p2}
Let $(\Lambda V,d)$ be a two-stage model of $X$, and $2\leq k_V\leq 3$. Then the torus rank of $X$ is less than ${\rm{dim}}~V^{\rm{odd}}-{\rm{dim}}~V^{\rm{even}}$.
\end{proposition}
\begin{proof}
Let $T^r$ acts almost freely on $X$, then we have a K-S extension
$$(\Lambda(x_1,x_2,\cdots,x_r),0)\rightarrow (\Lambda(x_1,x_2,\cdots,x_r)\otimes \Lambda V,D)\rightarrow (\Lambda V, d),$$
where $(\Lambda(x_1,x_2,\cdots,x_r)\otimes \Lambda V,D)$ is an elliptic Sullivan algebra. By \cite[p435]{Y-S}, we have
$r\leq {\rm{dim}}~V^{\rm{odd}}-{\rm{dim}}~V^{\rm{even}}$.

If $r={\rm{dim}}~V^{\rm{odd}}-{\rm{dim}}~V^{\rm{even}}$, by \cite[Proposition 32.10]{Y-S}, we have
$$(\Lambda(x_1,x_2,\cdots,x_r)\otimes \Lambda V,D)$$
is isomorphic to a pure Sullivan algebra. As $2\leq k_V\leq 3$, $d_{\tau}$ is quadratic. Then there is a pure Sullivan algebra such that $d_{\tau}\neq 0$, a contradiction.
\end{proof}
\section{Proof of Theorem \ref{1}}
\begin{proof}[Proof of Theorem \ref{1}]
Let $(\Lambda V,d)$ be a two-stage model of $X$, and $T^r$ acts almost freely on $X$.
If $k_{V}\leq 1$, then by Proposition \ref{p1}, we have
$${\rm{dim}}~H(\Lambda V,d)\geq 2^{{\rm{dim}}~V^{\rm{odd}}-{\rm{dim}}~V^{\rm{even}}}.$$
By \cite[p435]{Y-S}, we have $r\leq {\rm{dim}}~V^{\rm{odd}}-{\rm{dim}}~V^{\rm{even}}$,
so we have $${\rm{dim}}~H(\Lambda V,d)\geq 2^r.$$

If $2\leq k_V\leq 3$, then by Proposition \ref{p1} and Theorem \ref{p3}, we have 
$${\rm{dim}}~H(\Lambda V,d)> 2^{{\rm{dim}}~V^{\rm{odd}}-{\rm{dim}}~V^{\rm{even}}-1}.$$
By Proposition \ref{p2}, we have $r\leq {\rm{dim}}~V^{\rm{odd}}-{\rm{dim}}~V^{\rm{even}}-1$,
so we have $${\rm{dim}}~H(\Lambda V,d)> 2^r.$$

The desired result follows.
\end{proof}


\begin{thebibliography}{10}
\bibitem{H-S}S. Halperin, {\em Lectures on minimal models}, Mem. Soc. Math. France, \textbf{230}, 1983.
\bibitem{C-V}C. Allday, V. Puppe, {\em Cohomological Methods in Transformation Groups}, Cambridge Studies in Advanced Mathematics \textbf{32}, Cambridge University Press, 1993.
\bibitem{A-C-P}C. Allday, V. Puppe, {\em Bounds on the Torus Rank}, Volume {\textbf{1217}} of Lecture Notes in Math., 1-10. Springer, Berlin, 1986.
\bibitem{Y-S}Y. Felix, S. Halperin, J.-C. Thomas, {\em Rational homotopy theory}, Graduate Texts in Mathematics \textbf{205}, Springer-Verlag, New York, 2001.
\bibitem{B-G}B. Jessup, G. Lupton, {\em Free torus actions and two-stage spaces}, Math. Proc. Cambridge. \textbf{137} (01), 191-207, 2004.
\end{thebibliography}
\end{document}